\input amstex
\documentstyle{amsppt}
\NoBlackBoxes
\input bull-ppt
\keyedby{bull259e/PAZ}

\define\bcho{\partial{\scr{C}_0}}   
\define\7{{\bold{H}}}
\define\C{{\bold{C}}}                    
\define\QQ{{\bold{Q}}}                   
\define\1{\scr D}                        

\define\tr{\operatorname{Tr}}            
\define\arccosh{\operatorname{arccosh}}  
\define\teich{{\bold{T}_{1,1}}}   

\define\RR{{\bold{R}}}                      
\define\hS{{\widehat{\scr{S}}}}      
\define\2{\scr C}                   
                    %
\define\M{\scr{M}}           
             %
\define\pl{\operatorname{pl}}              
\define\Gmu{G_{\mu}}

\topmatter
\cvol{26}
\cvolyear{1992}
\cmonth{Jan}
\cyear{1992}
\cvolno{1}
\cpgs{141-146}
\title Pleating coordinates for the Teichm\"{u}ller\\ 
space of a punctured
      torus\endtitle
\shorttitle{Pleating coordinates for the Teichm\"uller 
space}
\author Linda Keen and Caroline Series\endauthor
\address Mathematics Department,
CUNY, Lehman College, Bronx, New York 10468\endaddress
\address Mathematics Institute, Warwick
University, Coventry CV4 7AL,
United Kingdom\endaddress
\thanks Research partially supported by NSF grant 
DMS-8902881\endthanks
\date March 29, 1991 and, in revised form, June 5, 
1991\enddate
\abstract We construct new coordinates for the
Teichm\"uller space Teich of a punctured torus into 
$\bold{R}
\times\bold{R}^+$. The coordinates depend on the 
representation of
Teich as a space of marked Kleinian groups $G_\mu$ that 
depend
holomorphically on a parameter $\mu$ varying in a simply 
connected
domain in $\bold{C}$. They describe the geometry of the 
hyperbolic
manifold $\bold{H}^3/G_\mu$; they reflect exactly the 
visual
patterns one sees in the limit sets of the groups 
$G_\mu$; and they
are directly computable from the generators of 
$G_\mu$.\endabstract
\endtopmatter
 
\document
\heading{1.} Introduction\endheading
 
In this paper we announce results that will appear in 
\cite{4}, in which
we construct a new  embedding in $\RR \times \RR^+$ of 
the Teichm\"{u}ller
space $\teich$  of a punctured torus. The pullbacks of 
the natural
coordinates in $\RR \times \RR^+$ have simple geometric 
interpretations when
$\teich$ is realized as a subset of $\C$ using the Maskit 
embedding. In the
Maskit embedding, points of $\teich$ correspond to marked 
Kleinian groups
$\{\Gmu\}$ that depend holomorphically on a parameter 
$\mu$ that varies in a
simply connected domain $\M$ in $\C$. Figure 1 (see p.\ 
142), drawn by David Wright, shows the domain $\M$. The
`coordinate grid' in $\M$ is the preimage under our 
embedding of the
horizontal and vertical lines in $\RR \times \RR^+$.
 
\pageinsert
\vfill
\centerline{\smc Figure 1. \rm The Maskit embedding with 
pleating coordinates.}
\endinsert

%
Our coordinates, which we call {\it pleating 
coordinates\/},
 have the following virtues: they relate directly to
the geometry of the hyperbolic manifold $\7^3/\Gmu$; 
 they reflect exactly
the visual patterns one sees in the limit sets of the 
groups $\Gmu$; and
they are directly computable from the generators of $\Gmu$.
The
definition of the Maskit embedding prescribes that the 
regular set of each
group $\Gmu$ should contain precisely one invariant 
component $\Omega _0(\Gmu)$.
The Riemann  surface $\Omega _0(\Gmu)/\Gmu$ is a 
punctured torus, and this
defines the correspondence between $\Gmu$ and a point in 
$\M$.  Ideally, one
would like to relate the geometry of the surface $\Omega 
_0(\Gmu)/\Gmu$ directly
to the parameter $\mu$. This seems to be a very hard 
problem. However, it
{\it is\/} possible to determine the relationship of 
$\mu$ to the geometry
of the hyperbolic manifold $\7^3/\Gmu$. This is the basic 
idea of
\cite{4}.

The boundary  of the convex hull of the limit set of a 
Kleinian group
acting on $\7^3$ carries an intrinsic hyperbolic metric 
and is a union of
pleated surfaces in the sense of Thurston (see 
\cite{9,1}). (Roughly
speaking, a pleated surface is a hyperbolic surface in a 
hyperbolic
3-manifold that is bent or pleated along some geodesic 
lamination called its
{\it pleating locus\/}.) We use the term {\it pleating 
coordinates\/} because
our embedding reflects the geometry of this pleating.
 
Since there is a
natural bijective correspondence between the connected 
components of the
regular set of $\Gmu$  and those of its convex hull 
boundary, exactly one of
these boundary components, $\bcho$ say, is invariant 
under the action of
$\Gmu$.   The two quotients $\bcho/\Gmu=\hS_{\mu}$ and 
$\Omega _0(\Gmu)/\Gmu$
are topologically, but not conformally, the same (see 
\cite{1,3}).  Thus
the pleated surface $\hS_{\mu}$  is  topologically a 
punctured torus
(hereafter denoted by $S$). The parameters we use in our 
embedding
reflect the
geometry, not of $\Omega _0(\Gmu)/\Gmu$, but of 
 $\hS_{\mu}$. In the rest of this
paper, we shall see how the pleating coordinates describe 
$\hS_{\mu}$ in
terms of its pleating locus $\pl(\mu)$.
%
\heading{2.} Pleating rays\endheading
 
The  `vertical' lines of the coordinate grid in
Figure 1 are the locus of points in $\M$ along which 
$\pl(\mu)$ is a
particular fixed geodesic lamination on $S$. We call 
these lines {\it
pleating rays\/}. They may be thought of as `internal 
rays' in $\M$ since
they play a role analogous to that of the `external rays' 
of the
Mandelbrot set in the study of the dynamics of quadratic 
polynomials.

Define a {\it projective measured lamination\/} to be a 
geodesic lamination together
with a projective class of transverse measures. The set 
of projective measured
laminations on $S$  is naturally identified with 
$\widehat{\RR}= \RR \cup
\{\infty\}$ (see \cite{8}). Because the pleating locus 
always carries a
natural transverse measure, the {\it bending measure\/} 
 (see \cite{9,1}),
we obtain for each $\mu$ a  projective  measured 
lamination on $S,$ also denoted
$\pl(\mu)$.   For each $\lambda \in \RR$ ($\lambda \neq 
\infty$), there is
a unique nonempty pleating ray
$$\mathop{{\cal P}_{\lambda}} = \{ \mu \in \M \colon
\pl(\mu)=\lambda \}.$$
 We show that  the ray $\mathop{{\cal P}_{\lambda}}$ is 
asymptotic to
the line $\Re \mu = 2 \lambda$ as $\Im \mu \rightarrow 
\infty$.

As is well known  \cite{8}, the simple closed curves on $S$
correspond exactly to $\QQ \,\cup\, \{\infty\}$. More 
precisely,
for each rational $p/q$, there is a unique free homotopy 
class $[\gamma_{p/q}]
\in \pi_1(S)$ for which the corresponding geodesic  on 
the unpunctured torus
is in the $(p,q)$- homology class, and for which the 
geodesic $\gamma_{p/q}$  on
 the
punctured torus is simple. We call the sets $
\mathop{{\cal P}_{{p}/{q}}}$, {\it rational
pleating rays\/}.

The group $\Gmu$ is an embedding of $\pi_1(S)$ in
$SL(2,\C)$, and the free homotopy class of $\gamma_{p/q}$ 
corresponds to a
conjugacy class of $\Gmu$ under this embedding. Choose 
$g_{p/q}(\mu)$  in
this conjugacy class.
\thm{Theorem 1}
The rational ray $\mathop{{\cal P}_{{p}/{q}}}$ is the 
unique branch of the locus $\{ \mu \in \C
\colon \tr  g_{p/q}(\mu) >2 \}$ that is asymptotically 
vertical as $\Im \mu
\rightarrow \infty$. This branch contains no 
singularities, and it intersects
$\partial{\M }$ in a unique point. At this point, $\tr 
g_{p/q}(\mu)=2$.
\ethm
%
 
It is not hard to show that the rational ray 
$\mathop{{\cal P}_{{p}/{q}}}$ must be
contained in the locus where $ \tr  g_{p/q}(\mu)$ is 
real.  The asymptotic
behavior of the trace polynomials is a straightforward 
 consequence of the
trace identities, as described in more detail below. 
 What is much more
interesting is that the pleating ray is  {\it 
precisely\/} that branch of
the real locus defined above. The key point in proving 
this is to show that
along $\mathop{{\cal P}_{{p}/{q}}}$, the invariant 
component $\Omega _0(\Gmu)$ is a  
{\it circle chain\/},
that is, a union of overlapping circles that fit together 
in a manner
reflecting the continued fraction expansion of $p/q$. 
There are two main
points:
for sufficiently large $c>0$, the ray $\mathop{{\cal 
P}_{{p}/{q}}}$
intersects the line $\Im \mu = c$ in a unique point, and,
circle chains
persist under continuous deformations of the group $\Gmu$ 
along 
$\mathop{{\cal P}_{{p}/{q}}}$.
The endpoint of the pleating ray represents a cusp group in
which the element $g_{p/q}(\mu)$ has become parabolic. In 
\cite{2}, we
prove that there are no other cusp groups corresponding 
to $g_{p/q}(\mu)$ in
$\partial{\M}$.

The circle chain patterns are visually apparent in computer
pictures of the limit sets of groups on the pleating ray 
and near to
$\partial{\M}$.  When $\mu$ reaches the endpoint of the 
pleating ray, the
overlapping circles become tangent. Such chains of 
tangent circles were
discovered by David Wright in the course of a computer 
investigation of
$\partial{\M}$, and our interest in them was the starting 
point of the
present work.

Wright obtained his striking pictures of $\partial{\M}$  by
using an inductive procedure related to the continued 
fraction
expansion of $p/q$ to
canonically choose   a
particular word $W_{p/q} \in \Gmu$ in the
     conjugacy class of the image
   of $[\gamma(p/q)]$.
He computed  $\tr W_{p/q}$ as a polynomial in
$\mu$ by means of the trace identities; and, by using his 
enumeration scheme
to give a systematic choice of initial point, he used 
 Newton's method to
find, for each $p/q$, a particular root of the equation 
$\tr
W_{p/q}(\mu)=2$. These solutions form the boundary curve 
in Figure 1.
%
\heading{3.} Pleating length\endheading

As $\mu$ moves down each rational ray, the length
in $\7^3$ of the pleating locus $\pl(\mu)$ provides a 
natural parameter. This
parameter, however, is not continuous as $\mu$ moves 
across rays. In fact,
if  $\mu_n \in \mathop{{\cal P}_{p_n/q_n}}$ converges to 
$\mu \in 
\mathop{{\cal P}_{\lambda}}$, where
$\lambda$ is not rational,  the hyperbolic lengths of the 
pleating loci
$\pl(\mu_n)$  always approach infinity. However, it {\it 
is\/} possible to
define a global length parameter that is continuous as 
$\mu$ moves across
rays by making a special choice of transverse measure, 
which we call the
{\it pleating measure\/}, for each projective measured 
lamination $\mu$ on $S$. We define the
{\it pleating length $PL(\mu)$} of $\Gmu$ to be the 
length of $\pl(\mu)$
with respect to the pleating measure of  $\pl(\mu)$. On a 
rational ray $
\mathop{{\cal P}_{{p}/{q}}}$, 
the pleating length of $\Gmu$ turns out to be  the 
hyperbolic length of
$\gamma_{p/q}(\mu)$ divided by the intersection number of 
$\gamma_{p/q}$ with
the fixed curve $\gamma_{\infty}$. The pleating length 
gives a natural
parameterization of the pleating rays, and the horizontal 
lines of the
coordinate grid in Figure 1 are lines of constant 
pleating length.

To prove continuity properties of the pleating measure 
and pleating length,
we use the continuous dependence on $\mu$ of  the 
hyperbolic structure of
the convex hull boundary, the pleating locus, and the 
bending measure. These
facts are also needed in the  proof of  Theorem 
1. We
prove all of these results in a more general setting in 
\cite{3}.
%
\heading{4.} The coordinates\endheading

It is apparent from Figure 1 that the partial
foliation of $\M$ by the rational rays should extend to a 
foliation by the
real rays $\mathop{{\cal P}_{\lambda}},\; \lambda \in 
\RR$. To show that it does, we
characterize the irrational pleating rays as the real 
loci of a family of
holomorphic functions.  The {\it complex translation 
length\/} of a loxodromic
element $g \in SL(2,\C)$ is defined as $2 \arccosh (\tr 
g)/2$ (see
\cite{9}). It follows from Theorem 
1 that on the
rational ray $\mathop{{\cal P}_{{p}/{q}}}$, the 
polynomial $\tr g_{p/q}(\mu)$ is real valued 
and,
hence, that the complex translation length is real. Since $
\mathop{{\cal P}_{{p}/{q}}}$ is
connected, we can choose a well-defined branch of the 
complex translation
length of $g_{p/q}(\mu)$ by specifying that it be real on $
\mathop{{\cal P}_{{p}/{q}}}$. We show
that the family of functions $\{L_{p/q} =1/q \arccosh 
(\tr g_{p/q}(\mu))
\}_{p/q \in \QQ}$ is normal in $\M$ and that on 
$\mathop{{\cal P}_{{p}/{q}}}$,
 the function
$L_{p/q}(\mu)$ coincides with the pleating length 
$PL(\mu)$.  Taking limits
in   ${\cal O}(\M)$,  the space of analytic functions on 
$\M$  with the
topology of uniform convergence on compact subsets,  we 
prove
%
\thm{Theorem 2}
The family $\{L_{p/q}\} $ extends to
a family $\{L_{\lambda}\}_{\lambda \in \RR}$ of  complex 
analytic functions
defined on $\M$,  such  that the function from $\M$ to 
$\RR$, given by $\mu
\mapsto PL(\mu)$, and the function from  $\RR$  to 
${\cal{O}}(\M)$, given by
$\lambda \mapsto L_{\lambda}$, are both continuous and 
such that  the
function $L_{\lambda}$ is real valued on $\mathop{{\cal 
P}_{\lambda}}$.
\ethm
%

That
the real rays are a codimension one foliation of $\M$ 
follows from
\thm{Theorem 3}
The real pleating ray $\mathop{{\cal P}_{\lambda}}$
is a connected component of the real locus of 
 $L_{\lambda}$ in $\M$. This
component contains no singularities and is asymptotic to 
$\Re {\mu} = 2
\lambda$ as $\Im {\mu} \rightarrow \infty$.
 \ethm

Our main theorem is
\thm{Theorem 4}
 The map from $\M$ to $ \RR \times
\RR^+$ defined by $\mu \mapsto (\pl(\mu),
PL(\mu))$ is a homeomorphism onto
its image.
\ethm
%

The fact that the map described in Theorem 4
is surjective will be proved elsewhere.

In a future paper, we expect to use the methods described 
here to give a
complete description of $\partial{\M}$ and of the 
approach to $\partial{\M}$
along the internal rays.  In particular, we hope to give 
proofs of
McMullen's theorems \cite{6,7}, conjectured by Bers, 
 that the cusp
groups are dense in $\partial{\M}$ and that 
$\partial{\M}$ is a Jordan
curve.  Although our work here relates to the punctured 
torus, most of the
techniques we have developed apply more generally. We 
plan to extend our
analysis to any union of surfaces of finite topological 
type. David Wright
has already produced computer pictures of the analogous 
coordinatization for
the (one complex dimensional) Riley slice of Schottky 
space, and the discussion in 
\cite{4} goes over to that situation (see \cite{5}).
%

\heading {} Acknowledgments\endheading

We wish to express our thanks to a number of
people. David Wright introduced us to this problem and 
has generously
allowed us to use his computer pictures. Curt McMullen 
has graciously shared
his ideas and work with us. We also want to thank 
Jonathan Brezin, David
Epstein, Michael
Handel, Steve Kerckhoff, Paddy Patterson, and Bill 
Thurston for many helpful
conversations throughout the course of this work. 
 Finally, we would like to
acknowledge the support of the NSF, the SERC in the 
United Kingdom, the Danish
Technical University, and the IMS at SUNY at Stonybrook.
%

\enddocument